\documentclass[12pt]{article}
\usepackage[english]{babel}
\usepackage[cp1251]{inputenc}
\usepackage{amsmath,latexsym}
\usepackage[dvips]{graphicx}
\usepackage[mathcal]{euscript}
\begin{document}
\textbf{Periodic Gibbs Measures for Models with Uncountable Set of
Spin Values on a Cayley Tree}\\

\textbf{ U. A. Rozikov $\cdot$ F. H. Haydarov}\\

U.\ A.\ Rozikov\\ Institute of mathematics, Tashkent, Uzbekistan.\\
e-mail: rozikovu@yandex.ru\\

 F.\ H.\ Haydarov\\ National University of Uzbekistan,
Tashkent, Uzbekistan.\\
e-mail: haydarov\_ imc@mail.ru\\

 \textbf{Abstract}  We consider models with nearest-neighbor interactions and
with the set $[0,1]$ of spin values, on a Cayley tree of order
$k\geq 1$. We study periodic Gibbs measures of the model with
period two. For $k=1$ we show that there is no any periodic Gibbs
measure. In case $k\geq 2$ we get a sufficient condition on
Hamiltonian of the model with uncountable set of spin values under
which the model have not any periodic Gibbs measure with period
two. We construct several models which have at least two periodic
Gibbs measures.\\


\textbf{1 Introduction} \\

 The structure of the lattice (graph) plays an
important role in investigations of spin systems. For example, in
order to study the phase transition problem for a system on $Z^d$
and on Cayley tree there are two different methods: Pirogov-Sinai
theory on $Z^d$, Markov random field theory and recurrent
equations of this theory on Cayley tree. In [1-5,8,11-13, 15-17]
for several models on Cayley tree, using the Markov random field
theory Gibbs measures are described.

These papers are devoted to models with a $ finite $ set of spin
values. Mainly were shown that these models have finitely many
translation-invariant  and uncountable numbers  of the
non-translation-invariant extreme Gibbs measures. Also for several
models (see, for example, [6,8,12]) it were proved that there
exist three periodic  Gibbs  measures (which are invariant with
respect to normal  subgroups  of  finite index of the group
representation of the Cayley tree) and there are uncountable
number of non-periodic Gibbs measures.

In [7] the Potts model with a countable set of spin values on a
Cayley tree is considered and it was showed that the set of
translation-invariant splitting Gibbs measures of the model
contains at most one point, independently on parameters of the
Potts model with countable set of spin values on the Cayley tree.
This is a crucial difference from the models with a finite set of
spin values, since the last ones may have more than one
translation-invariant Gibbs measures.

In [3], [4], [13] models with an uncountable set of spin values
are considered. Our paper is continuation of these papers.

The paper is organized as follows. Section 2 introduces the main
definitions. In Sect.3 we prove non-existance Gibbs measures with
period two on Cayley tree of order one. In Sect.4. the
Hammerstein's nonlinear integral equation is presented. In Sect.5.
we give a sufficient condition on Hamiltonian of the model have
not any periodic Gibbs measure. In Sect 6,7 and 8 the existance of
at least two periodic Gibbs measures for several models with
uncountable set of spin values are proved respectively in cases
$k=2$, $k=3$, $k\geq 4.$ In Sect.9. the existance of at least four
periodic Gibbs measures for the models with uncountable set of
spin values are proved in cases $k\geq k_{0}.$ \\

\textbf{2 Preliminaries }\\

A Cayley tree $\Gamma^k=(V,L)$ of order $k\geq 1$ is an infinite
homogeneous tree, i.e., a graph without cycles, with exactly $k+1$
edges incident to each vertices. Here $V$ is the set of vertices
and $L$ that of edges (arcs).

Consider models where the spin takes values in the set $[0,1]$,
and is assigned to the vertexes of the tree. For $A\subset V$ a
configuration $\sigma_A$ on $A$ is an arbitrary function
$\sigma_A:A\to [0,1]$. Denote $\Omega_A=[0,1]^A$ the set of all
configurations on $A$. A configuration $\sigma$ on $V$ is then
defined as a function $x\in V\mapsto\sigma (x)\in [0,1]$; the set
of all configurations is $[0,1]^V$. The (formal) Hamiltonian of
the model is :

$$ H(\sigma)=-J\sum_{\langle x,y\rangle\in L}
\xi_{\sigma(x),\sigma(y)}, \eqno(2.1)$$
 where $J \in R\setminus \{0\}$ and
$\xi: (u,v)\in [0,1]^2\to \xi_{u,v}\in R$ is a given bounded,
measurable function. As usually, $\langle x,y\rangle$ stands for
nearest neighbor vertices.

Let $\lambda$ be the Lebesgue measure on $[0,1]$.  On the set of
all configurations on $A$ the a priori measure $\lambda_A$ is
introduced as the $|A|$ fold product of the measure $\lambda $.
Here and further on $|A|$ denotes the cardinality of $A$.   We
consider a standard sigma-algebra ${\mathcal B}$ of subsets of
$\Omega=[0,1]^V$ generated by the measurable cylinder subsets.
 A probability measure $\mu$ on $(\Omega,{\mathcal B})$
is called a Gibbs measure (with Hamiltonian $H$) if it satisfies
the DLR equation, namely for any $n=1,2,\ldots$ and
$\sigma_n\in\Omega_{V_n}$:
$$\mu\left(\left\{\sigma\in\Omega :\;
\sigma\big|_{V_n}=\sigma_n\right\}\right)= \int_{\Omega}\mu ({\rm
d}\omega)\nu^{V_n}_{\omega|_{W_{n+1}}} (\sigma_n),$$ where
$\nu^{V_n}_{\omega|_{W_{n+1}}}$ is the conditional Gibbs density
$$ \nu^{V_n}_{\omega|_{W_{n+1}}}(\sigma_n)=\frac{1}{Z_n\left(
\omega\big|_{W_{n+1}}\right)}\exp\;\left(-\beta H
\left(\sigma_n\,||\,\omega\big|_{W_{n+1}}\right)\right),
$$
and $\beta={1\over T}$, $T>0 $ is temperature. Here and below,
$W_l$ stands for a `sphere' and $V_l$ for a `ball' on the tree, of
radius $l=1,2,\ldots$, centered at a fixed vertex $x^0$ (an
origin):
$$W_l=\{x\in V: d(x,x^0)=l\},\;\;V_l=\{x\in V: d(x,x^0)\leq l\};$$
and
$$L_n=\{\langle x,y\rangle\in L: x,y\in V_n\};$$
distance $d(x,y)$, $x,y\in V$, is the length of (i.e. the number
of edges in) the shortest path connecting $x$ with $y$.
$\Omega_{V_n}$ is the set of configurations in $V_n$ (and
$\Omega_{W_n}$ that in $W_n$; see below). Furthermore,
$\sigma\big|_{V_n}$ and $\omega\big|_{W_{n+1}}$ denote the
restrictions of configurations $\sigma,\omega\in\Omega$ to $V_n$
and $W_{n+1}$, respectively. Next, $\sigma_n:\;x\in V_n\mapsto
\sigma_n(x)$ is a configuration in $V_n$ and
$H\left(\sigma_n\,||\,\omega\big|_{W_{n+1}}\right)$ is defined as
the sum $H\left(\sigma_n\right)+U\left(\sigma_n,
\omega\big|_{W_{n+1}}\right)$ where
$$H\left(\sigma_n\right)
=-J\sum_{\langle x,y\rangle\in
L_n}\xi_{\sigma_n(x),\sigma_n(y)},$$
$$U\left(\sigma_n,
\omega\big|_{W_{n+1}}\right)= -J\sum_{\langle x,y\rangle:\;x\in
V_n, y\in W_{n+1}} \xi_{\sigma_n(x),\omega (y)}.$$ Finally,
$Z_n\left(\omega\big|_{W_{n+1}}\right)$ stands for the partition
function in $V_n$, with the boundary condition
$\omega\big|_{W_{n+1}}$:
$$Z_n\left(\omega\big|_{W_{n+1}}\right)=
\int_{\Omega_{V_n}} \exp\;\left(-\beta H
\left({\widetilde\sigma}_n\,||\,\omega
\big|_{W_{n+1}}\right)\right)\lambda_{V_n}(d{\widetilde\sigma}_n).$$

Due to the nearest-neighbor character of the interaction, the
Gibbs measure possesses a natural Markov property: for given a
configuration $\omega_n$ on $W_n$, random configurations in
$V_{n-1}$ (i.e., `inside' $W_n$) and in $V\setminus V_{n+1}$
(i.e., `outside' $W_n$) are conditionally independent.

We use a standard definition of a periodic measure (see,
e.g.[9],[15]).
 The main object of study in this
paper are periodic Gibbs measures for the model (2.1) on Cayley
tree. In $[13]$ the problem of description of such measures was
reduced to the description of the solutions of a nonlinear
integral equation. For finite and countable sets of spin values
this argument is well known (see, e.g. [1-7, 9-10, 15-17]).

Write $x<y$ if the path from $x^0$ to $y$ goes through $x$. Call
vertex $y$ a direct successor of $x$ if $y>x$ and $x,y$ are
nearest neighbors. Denote by $S(x)$ the set of direct successors
of $x$. Observe that any vertex $x\ne x^0$ has $k$ direct
successors and $x^0$ has $k+1$.

Let $h:\;x\in V\mapsto h_x=(h_{t,x}, t\in [0,1]) \in R^{[0,1]}$ be
mapping of $x\in V\setminus \{x^0\}$.  Given $n=1,2,\ldots$,
consider the probability distribution $\mu^{(n)}$ on
$\Omega_{V_n}$ defined by

$$\mu^{(n)}(\sigma_n)=Z_n^{-1}\exp\left(-\beta H(\sigma_n)
+\sum_{x\in W_n}h_{\sigma(x),x}\right). \eqno(2.2)$$
 Here, as before, $\sigma_n:x\in V_n\mapsto
\sigma(x)$ and $Z_n$ is the corresponding partition function:

$$Z_n=\int_{\Omega_{V_n}} \exp\left(-\beta H({\widetilde\sigma}_n)
+\sum_{x\in W_n}h_{{\widetilde\sigma}(x),x}\right)
\lambda_{V_n}(\widetilde{\sigma}_{n}). \eqno(2.3)$$

The probability distributions $\mu^{(n)}$ are compatible if for
any $n\geq 1$ and $\sigma_{n-1}\in\Omega_{V_{n-1}}$:

$$\int_{\Omega_{W_n}}\mu^{(n)}(\sigma_{n-1}\vee\omega_n)\lambda_{W_n}(d(\omega_n))=
\mu^{(n-1)}(\sigma_{n-1}). \eqno(2.4)$$
 Here
$\sigma_{n-1}\vee\omega_n\in\Omega_{V_n}$ is the concatenation of
$\sigma_{n-1}$ and $\omega_n$. In this case there exists a unique
measure $\mu$ on $\Omega_V$ such that, for any $n$ and
$\sigma_n\in\Omega_{V_n}$, $\mu \left(\left\{\sigma
\Big|_{V_n}=\sigma_n\right\}\right)=\mu^{(n)}(\sigma_n)$.\\

\textbf{Definition 2.1} The measure $\mu$ is called {\it splitting
Gibbs measure} corresponding to Hamiltonian (2.1) and function
$x\mapsto h_x$, $x\neq x^0$.

 The following statement describes conditions on $h_x$ guaranteeing compatibility
of the corresponding distributions $\mu^{(n)}(\sigma_n).$\\

 \textbf{Proposition 2.2}[13] The probability distributions
$\mu^{(n)}(\sigma_n)$, $n=1,2,\ldots$, in (2.2) {\sl are
compatible iff for any $x\in V\setminus\{x^0\}$ the following
equation holds:

$$ f(t,x)=\prod_{y\in S(x)}{\int_0^1\exp(J\beta\xi_{t,u})f(u,y)du
\over \int_0^1\exp(J\beta{\xi_{0,u}})f(u,y)du}. \eqno(2.5)$$
 Here,
and below $f(t,x)=\exp(h_{t,x}-h_{0,x}), \ t\in [0,1]$ and
$du=\lambda(du)$ is the Lebesgue measure.}\\

From Proposition 2.2 it follows that for any $h=\{h_x\in
R^{[0,1]},\ \ x\in V\}$ satisfying (2.5) there exists a unique
Gibbs measure $\mu$ and vice versa. However, the analysis of
solutions to (2.5) is not easy. This difficulty depends on the
given function $\xi$.\\

$x\in V$ is called even (odd) if $d(x,x_{0})$ - even (odd). In
this paper we shall study special periodic solutions to (2.5),
which are in the form $f(t,x)=f(t)$ if $x$ - even and
$f(t,x)=g(t)$ if $x$ - odd. For such functions equation (2.5) can
be written as

$$f(t)=\left({\int_0^1K(t,u)g(u)du\over \int_0^1
K(0,u)g(u)du}\right)^k, \,\  g(t)=\left({\int_0^1K(t,u)f(u)du\over
\int_0^1 K(0,u)f(u)du}\right)^k, \eqno(2.6)$$
 where $K(t,u)=\exp(J\beta \xi_{tu}),
f(t),g(t)>0, t,u\in [0,1].$

We put
$$C^+[0,1]=\{ \varphi \in C[0,1]: \varphi(x)\geq 0\}.$$
We are interested to positive continuous solutions to (2.6), i.e.
such that

$$f,g \in C_0^+[0,1]=\{\varphi\in C[0,1]: \varphi(x)\geq
0\}\setminus \{\theta\equiv 0\}.$$
 Define the operator $A_{k}:C^{+}_{0}[0,1]\rightarrow C^{+}_{0}[0,1]$ by
$$(A_{k}f)(t)=\left[\frac{(Wf)(t)}{(Wf)(0)}\right]^{k}, \,\
k\in\mathbf{N},$$\\
where $W:C[0,1]\rightarrow C[0,1]$ is linear operator, which is
defined by :
 $$(Wf)(t)=\int^{1}_{0}K(t,u)f(u)du, \eqno(2.7)$$
and defined the linear functional $\omega: C[0,1]\rightarrow R$ by
$$\omega(f)=(Wf)(0)=\int_{0}^{1}K(0,u)f(u)du.$$
 Then Eq.(2.6) can be written as

$$A_{k}f=g, \,\  A_{k}g=f, \,\,\,\,\,\,\,\,\,\,\,\,\,\,\ f,g\in C^{+}_{0}[0,1] \eqno(2.8)$$\\
\textbf{3 Non Existence of periodic Gibbs Measures for the Model (2.1): Case $k=1$}.\\

 At first we are going to consider (2.8) for $k=1$. The system of equations (2.8) is equivalent
 to linear equations

  $$(Wf)(t)=w(f)g(t), \,\,\ (Wg)(t)=w(g)f(t) , \,\,\,\,\,\  f,g\in
  C^{+}[0,1] \eqno(3.1)$$

 \textbf{Lemma 3.1} Let $(f,g)$ satisfies (3.1) with $f\neq g$, $\delta_{0}=sup \{\delta\in (0,\infty):f-\delta
 g>0\}.$\\
Then $W(f-\delta_{0}g)>0.$

\emph{Proof} We have $ f-\delta_{0}g\geq 0$ \,\ $\Rightarrow$ $
W(f-\delta_{0}g)\geq 0 .$ \,\ Suppose $W(f-\delta_{0}g)=0$ then
 $$f-\delta_{0}g\equiv 0 \,\,\ \Rightarrow \frac{f(t)}{g(t)}=\delta_{0}, \,\,\,\  t\in [0,1].$$
For \,\,\ $t=0$
 $$g(0)=\frac{(Wf)(0)}{w(f)}=1=\frac{(Wg)(0)}{w(g)}=f(0).$$
Then $\delta_{0}=1.$ \,\ This contradicts to $f\neq g.$ \,\ Thus
we have proved $W(f-\delta_{0}g)>0.$ \\

 \textbf{Theorem 3.2} \,\,\ The system of equations $(A_{1}f)(t)=g(t)$ ,
 $(A_{1}g)(t)=f(t)$ has not any solution $(f,g)\in
 (C^{+}[0,1])^{2}$ with $f\neq g.$\\

  \emph{Proof} \,\,\ Let $(f_{1}(x),g_{1}(x))$ be a solution of the system of
  equations:
   $$(A_{1}f)(t)=g(t), \,\,\  (A_{1}g)(t)=f(t).$$
Then  $$(Wf_{1})(t)=w(f_{1})g_{1}(t)
  ,\,\
 (Wg_{1})(t)=w(g_{1})f_{1}(t).$$
Put $$ \lambda_{1}=w(f_{1}), \,\,\  \lambda_{2}=w(g_{1}) \,\,\,\
\Rightarrow \,\ \lambda_{i}>0 , \,\  i\in\{1,2\}.$$ Denote
$$\delta_{1}=sup \{\delta\in (0,\infty): f-\delta g>0\},\,\,\
\delta_{2}=sup \{\delta\in (0,\infty): g-\delta f>0\}.$$
By Lemma
3.1
  $$\lambda_{1}g(t)-\lambda_{2}\delta_{1}f(t)=W(f-\delta_{1} g)>0,$$
  $$\lambda_{2}f(t)-\lambda_{1}\delta_{2}g(t)=W(g-\delta_{2} f)>0.$$
  Hence $$\frac{\lambda_{2}}{\lambda_{1}}\delta_{1}<\frac{g(t)}{f(t)} , \,\,\
   \frac{\lambda_{1}}{\lambda_{2}}\delta_{2}<\frac{f(t)}{g(t)},  \,\,\ t\in[0,1].$$
   There exists $t_{0} , t_{1}\in[0,1]$ such that $\delta_{2}=\frac{g(t_{0})}{f(t_{0})}$
and $\delta_{1}=\frac{f(t_{1})}{g(t_{1})}.$ Then
   $$\frac{g(t)}{f(t)}\geq\frac{g(t_{0})}{f(t_{0})}=\delta_{2}>\frac{\lambda_{2}}{\lambda_{1}}\delta_{1}
   , \,\,\,\ \frac{f(t)}{g(t)}\geq\frac{f(t_{1})}{g(t_{1})}=\delta_{1}>\frac{\lambda_{1}}{\lambda_{2}}\delta_{2}. $$
  Thus we have $\frac{\lambda_{2}}{\lambda_{1}}\delta_{1}<\delta_{2}$  and
 $\frac{\lambda_{1}}{\lambda_{2}}\delta_{2}<\delta_{1}$
 \,\,\  this is a contradiction. \\

 \textbf{4  The Hammerstain's nonlinear equation.}\\

For every $k\in N$ we consider an integral operator $H_{k}$ acting
 in $C^{+}[0,1]$ as
  $$ (H_{k}f)(t)=\int^{1}_{0}K(t,u)f^{k}(u)du \eqno(4.1)$$
 If $k\geq2$ then the operator $H_{k}$ is a nonlinear operator which is
 called Hammerstain's operator of order $k$.
  For a nonlinear homogenous operator A  it is known that if there
  are positive solutions of the operator A then the number of the
  positive solutions are continuum. (see[10] , p.186).

Denote
$$\mathcal M_{0}=\{\varphi\in C^{+}[0,1]:\varphi(0)=1\}.$$

\textbf{Lemma 4.1} The system of equations:
 $$(A_{k}f)(t)=g(t), \,\,\,\,\  (A_{k}g)(t)=f(t) \,\,\,\,\,\,\,\  k\geq2. \eqno(4.2)$$
  has a positive solution iff the system of equations:
 $$ (H_{k}f)(t)=\lambda_{1}g(t), \,\,\,\,\ (H_{k}g)(t)=\lambda_{2}f(t), \,\,\,\,\,\ k\geq2 \eqno(4.3)$$
 has a positive solution in $(\mathcal M_{0})^{2}.$\\

 \emph{Proof}  {\sl Necessity.} Let $(f_{0} , g_{0}) \in (C^{+}_{0}[0,1])^{2}$
 be a solution of the system of equations (4.2). We have
 $$ (Wf_{0})(t)=w(f_{0})\sqrt[k]{g_{0}(t)}, \,\,\,\ (Wg_{0})(t)=w(g_{0})\sqrt[k]{f_{0}(t)}.$$
 From this equality we get
$$ (H_{k}f_{1})(t)=\lambda_{1}g_{1}(t), \,\,\,\,\ (H_{k}g_{1})(t)=\lambda_{2}f_{1}(t).$$
where $f_{1}(t)=\sqrt[k]{f_{0}(t)}$ ,
$g_{1}(t)=\sqrt[k]{g_{0}(t)}$ \,\ and  \,\
$\lambda_{1}=w(f_{0})>0$ , $\lambda_{2}=w(g_{0})>0.$ It is easy to
see that $(f_{1},g_{1})\in (\mathcal M_{0})^{2}$

{\sl Sufficiency.} Let $k\geq2$ and $(f_{1} , g_{1}) \in (\mathcal
M_{0})^{2}$ be a solution of the system (4.3). From $f_{1}(0)=1$ ,
$g_{1}(0)=1$ we get
$$ 1=g_{1}(0)=(H_{k}f_{1})(0)=w(f_{1}^{k}), \,\,\,\,\ 1=f_{1}(0)=(H_{k}g_{1})(0)=w(g_{1}^{k}).$$
Then
 $$f_{1}=\frac{H_{k}g_{1}}{w(g_{1}^{k})}, \,\,\,\,\  g_{1}=\frac{H_{k}f_{1}}{w(f_{1}^{k})}. $$
 From this equalities we get $A_{k}f_{0}=g_{0},$ $A_{k}g_{0}=f_{0}$
 with $f_{0}=f_{1}^{k}\in C^{+}_{0}[0,1],$ $g_{0}=g_{1}^{k}\in
 C^{+}_{0}[0,1].$ This completes the proof.\\

\textbf{Lemma 4.2} The system of equations (4.3)
  has a positive solution iff the system of equations:
 $$ (H_{k}f)(t)=g(t), \,\,\,\,\ (H_{k}g)(t)=f(t), \,\,\,\,\,\ k\geq2 \eqno(4.4)$$
 has a positive solution. \\

 \emph{Proof} {\sl Necessity.} Let $(f_{0}(t),g_{0}(t))$ be
 a positive solution of the system (4.3). Define functions:
 $$ f_{1}(t)=\frac{1}{C_{1}}f_{0}(t), \,\,\,\,\  g_{1}(t)=\frac{1}{C_{2}}g_{0}(t).$$
Then
 $$(H_{k}f_{1})(t)=\frac{1}{C_{1}^{k}}(H_{k}f_{0})(t)=\frac{\lambda_{1}}{C_{1}^{k}}g_{0}(t)=
 \frac{\lambda_{1}C_{2}}{C_{1}^{k}}g_{1}(t).$$
 Put
 $$C_{1}=(\lambda_{1})^{\frac{1}{k+1}}(\lambda_{1}\lambda_{2})^{\frac{1}{k^{2}-1}},  \,\,\,\,\,\
 C_{2}=(\lambda_{2})^{\frac{1}{k+1}}(\lambda_{1}\lambda_{2})^{\frac{1}{k^{2}-1}}. \eqno(4.5)$$
We have
 $$(H_{k}f_{1})(t)=\frac{\lambda_{1}(\lambda_{2})^{\frac{1}{k+1}}(\lambda_{1}\lambda_{2})^{\frac{1}{k^{2}-1}}}
 {\lambda_{1}^{\frac{k}{k+1}}(\lambda_{1}\lambda_{2})^{\frac{k}{k^{2}-1}}}g_{1}(t)=
 \frac{(\lambda_{1}\lambda_{2})^{\frac{1}{k+1}}}{(\lambda_{1}\lambda_{2})^{\frac{k-1}{k^{2}-1}}}g_{1}(t)=g_{1}(t).$$
Similarly we get
$$H_{k}g_{1}(t)=f_{1}(t).$$
{\sl Sufficiency.} Let $(f_{1}(t),g_{1}(t))$ be a positive
solution of the system (4.4). We get
$$f_{0}(t)=C_{1}f_{1}(t), \,\,\,\ g_{0}(t)=C_{2}g_{1}(t).$$
where $C_{1},\,\ C_{2}$ are given in (4.5)
 It is easy to verify
$$(H_{k}f_{0})(t)=\lambda_{1}g_{0}(t), \,\,\,\ (H_{k}g_{0})(t)=\lambda_{2}f_{0}(t), \,\,\,\,\,\,\,\,\ k\geq2. $$
This completes the proof. \\

 Denote
$$\mathcal K=\left\{f\in C^+[0,1]: M \cdot\min_{t\in [0,1]}f(t)\geq m\cdot\max_{t\in [0,1]}f(t)\right\},$$
$$\mathcal P_k=\left\{\varphi\in C[0,1]: {m\over M} \cdot \left(1\over M\right)^{1\over k-1}\leq \varphi(t)
\leq {M\over m}\cdot\left(1\over m\right)^{1\over k-1}\right\}, \,
k\geq 2, $$ where
$$M=\max_{t,u\in [0,1]^{2}}K(t,u) , \,\,\,\,\  m=\min_{t,u\in [0,1]^{2}}K(t,u).$$ \\

 \textbf{Proposition 4.3} Let $k\geq 2$. Then \\

a) $H_k(C^+[0,1])\subset \mathcal K.$\\

b) If $(f_0, g_{0})\in (C_0^+[0,1])^{2}$ is a solution of the
system (4.4) then $(f_0,  g_{0})\in (\mathcal P_k)^{2}$.\\

 \emph{Proof} \,\ a) Let $\varphi \in H_k(C^+[0,1])$ be an arbitrary function. There exists a function $\psi\in C^+[0,1]$
such that $\varphi=H_k\psi$. Since $\varphi$ is continuous on
$[0,1]$, there are $t_1,t_2\in [0,1]$ such that
$$\varphi_{\min}=\min_{t\in[0,1]}\varphi(t)=\varphi(t_1)=(H_k\psi)(t_1),$$
$$\varphi_{\max}=\max_{t\in[0,1]}\varphi(t)=\varphi(t_2)=(H_k\psi)(t_2).$$
Hence
$$\varphi_{\min}\geq m\int^1_0\psi^k(u)du\geq m\int^1_0{K(t_2,u)\over M}\psi^k(u)du={m\over M}\varphi_{max},$$
i.e. $\varphi \in \mathcal K$.

 b) Let $(f,g)\in (C_0^+[0,1])^{2}$ be a solution of the system (4.4). Then we have
 $$\|f\|\leq M\|g\|^{k}, \,\,\,\ \|g\|\leq
M\|f\|^k,$$ where
$$||f||=\max_{t\in[0,1]}|f(t)|.$$
Hence
$$\|f\|\leq M^{k+1}\|f\|^{k^{2}}.$$ Consequently
$$\|f\|\geq \left(1\over M\right)^{1\over k-1}.$$
By the property a)
$$f(t)\geq f_{\min}=\min_{t\in [0,1]}f(t)\geq {m\over M}\|f\|.$$
Consequently
$$f(t)\geq {m\over M}\left({1\over M}\right)^{1\over k-1}.$$
Also we have
$$f(t)=(H_kg)(t)\geq m\int^1_0g^k(u)du\geq m g_{\min}^k$$
and
$$g(t)=(H_kf)(t)\geq m\int^1_0f^k(u)du\geq m f_{\min}^k.$$
Then $f_{\min}\geq m g_{\min}^k$ ,  $g_{\min}\geq m f_{\min}^k$
 $\Rightarrow f_{\min}\geq m^{k+1} f_{\min}^{k^{2}}$
i.e.
$$f_{\min}\leq \left({1\over m}\right)^{1\over k-1}.$$
 By the property a)
$$f(t)\leq f_{\max}\leq {M\over m}f_{\min}\leq {M\over m}\left({1\over m}\right)^{1\over k-1}.$$ Thus we have
$f\in \mathcal P_k$. Similarly one can prove that $g\in \mathcal
P_k.$  \\

 \textbf{Lemma 4.4} Let $f\neq g$. Put
 $$\delta_{1}=sup\{\delta\in [0,\infty): f(t)-\delta g(t)\in C^+[0,1]\}$$
 and
 $$\delta_{2}=sup\{\delta\in [0,\infty): g(t)-\delta f(t)\in C^+[0,1]\}.$$
If max$\{\delta_{1},\delta_{2}\}\geq1$ then $(f,g)\in
(C_0^+[0,1])^{2}$ can not be solution to the system (4.4). \\

 \emph{Proof} Let $(f_{1},g_{1})\in(C_{0}^{+}[0,1])^{2}$ be a solution of system (4.4) and
 assume max$\{\delta_{1},\delta_{2}\}=\delta_{1}\geq1$ (the case max$\{\delta_{1},\delta_{2}\}=\delta_{2}$ is
 similar).
 Then
 $$g(t)-\delta_{1}^{k}f(t)=\int^{1}_{0}K(x,t)(f^{k}(x)-\delta_{1}^{k}g^{k}(x))dx\geq0.$$
 There exists  $t\exists_{0}\in[0,1]$ such that
 $\delta_{1}=\frac{g(t_{0})}{f(t_{0})}.$
 Moreover we have
  $$\frac{g(t)}{f(t)}\geq\delta_{1}^{k}, \,\,\,\,\,\,\  {t\in[0,1].}$$
Then
 $$\delta_{1}=\frac{g(t_{0})}{f(t_{0})}\geq\delta_{1}^{k} \Rightarrow \delta_{1}=1.$$
 It is clear that if $\delta_{1}=1$ then $f(x)=g(x).$  But this contradicts to $f(x)\neq g(x).$  \\

  \textbf{Theorem 4.5} Let $(f_{1}(t),g_{1}(t))$ be a solution of
  system (4.4) with $f_{1}\neq g_{1}$.
  Put $\varphi(t)=f_{1}(t)-g_{1}(t)$. The function
  $\varphi(t)$  changes its sign in $[0,1]$.\\

\emph{Proof}  Assume that $f_{1}(t)-g_{1}(t)\geq0$ (the case
$g_{1}(t)-f_{1}(t)\geq0$ is similar).\\
Consider
 $$\varphi_{\delta}^{(1)}(t)=f_{1}(t)-\delta g_{1}(t), \,\,\ \varphi_{\delta}^{(2)}(t)=g_{1}(t)-\delta
 f_{1}(t),\,\,\,\,\ \delta\in [0,\infty).$$
 Put
$$\delta_{1}=sup\{\delta\in [0,\infty):\varphi_{\delta}^{(1)}(t)\in C^+[0,1]\}$$
 and
 $$\delta_{2}=sup\{\delta\in [0,\infty): \varphi_{\delta}^{(2)}(t)\in C^+[0,1]\}.$$
 One can easily check that $\delta_{1}\geq1.$
 By Lemma 4.4 $(f_{1}(t) , g_{1}(t))$ can not be solution of
 system (4.4). This contradicts our assumption $f_{1}(t)-g_{1}(t)\geq0.$ \\

\textbf{5 Non Existence of periodic Gibbs Measures for
Model(2.1): Case $k\geq2$.}\\

\textbf{Lemma 5.1} Assume function $\varphi\in C[0,1]$ changes its
sign on $[0,1]$. Then for every $a\in R$ the following inequality
holds
$$\|\varphi_a\|\geq {1\over n+1}\|\varphi\|, \ \ n\in \mathbf{N},$$
where $\varphi_a=\varphi_a(t)=\varphi(t)-a, \, t\in [0,1].$ (see [3]. p.9)\\

\textbf{Proposition 5.2} Let $k\geq 2$. If the kernel $K(t,u)$
satisfies the condition $$ \left(M\over m\right)^k-\left(m\over
M\right)^k<{1\over k}, \eqno(5.1)$$ then the system (4.4) has
not any solution $(f,g)$ in $(C_0^+[0,1])^{2}$ with $f\neq g.$\\

 \emph{Proof}
 Assume that there is a solution
$(f_1 , g_1)\in (C_0^+[0,1])^{2}$. Denote $h(t)=f_1(t)-g_1(t)$.
Then by Theorem 4.5 the function $h(t)$ changes its sign on
$[0,1]$. By Lemma 5.1  we get
$$\max_{t\in[0,1]}\left|h(t)+{k\over 2}(\gamma_1+\gamma_2)\int^1_0f(s)ds\right|\geq{1\over 2}\|h\|,$$
where
$$\gamma_1=\left(m\over M\right)^k, \ \ \gamma_2=\left(M\over m\right)^k.$$
By a mean value Theorem we have
$$-h(t)=\int^1_0K(t,u)k\xi^{k-1}(u)h(u)du,$$
here $\xi\in C^+[0,1]$ and
$$\min\{f_1(t), g_1(t)\}\leq \xi(t)\leq \max\{f_1(t), g_1(t)\}, \, t\in [0,1].$$
By Proposition 4.3 we have $\xi\in \mathcal P_k$, i.e.
$${m\over M}\left(1\over M\right)^{1\over k-1}\leq \xi(t)\leq {M\over m}\left(1\over m\right)^{1\over k-1},
 \, t\in [0,1].$$
Hence
$$\gamma_1\leq K(t,u)\xi^{k-1}(u)\leq \gamma_2,\, t,u\in[0,1].$$
Therefore
$$\left|k\cdot K(t,u)\xi^{k-1}(u)-k{\gamma_1+\gamma_2\over 2}\right|\leq k{\gamma_2-\gamma_1\over 2}.$$
Then $$ \left|h(t)-{k\over
2}(\gamma_1+\gamma_2)\int^1_0h(u)du\right|\leq {k\over 2}
(\gamma_2-\gamma_1)\|h\|. \eqno(5.2)$$
 Assume the kernel $K(t,u)$ satisfies
the condition (5.1). Then $k(\gamma_2-\gamma_1)<1$ and the
inequality (5.2) contradicts to Lemma 5.1.
This completes the proof. \\

 \textbf{Proposition 5.3} Let $k\geq 2$. Let the kernel $K(t,u)$ satisfies
the condition (5.1). For every $\lambda_{1}>0$,
 $\lambda_{2}>0$ the Hammerstein's system of equations
 $$H_kf=\lambda_{1}g , \,\,\,\ H_kg=\lambda_{2}f \eqno(5.3)$$
 has not solution $(f,g)\in (C_0^+[0,1])^{2}, f\neq g.$\\

 \emph{Proof}  By Lemma 4.2 the system of equations (5.3)
is equivalent to the following system of equations
$$\int^1_0K(t,u)f_{1}^k(u)du=g_{1}(t), \,\,\,\ \int^1_0K(t,u)g_{1}^k(u)du=f_{1}(t) \eqno(5.4)$$
 By Theorem 5.1 the equation (5.4) has not solution in
$(C_0^+[0,1])^{2}$.
 Hence the equation (5.3) has not solution in $(C_0^+[0,1])^{2}$.\\

\textbf{Theorem 5.4} Let $k\geq 2$. If the kernel $K(t,u)$
satisfies the condition (5.1), then the system of equations (2.6)
has not solution in $(C_0^+[0,1])^{2}, f\neq g.$\\

\emph{Proof} Assume there is solution $(f_1,g_1)\in
(C^+[0,1])^{2}$,
 i.e.
 $$A_kf_{1}=g_{1}, \,\,\ A_{k}g_{1}=f_{1}.$$
 By Lemma 4.1 the functions $f_{2}(t)=\sqrt[k]{f_1(t)}$ and $g_{2}(t)=\sqrt[k]{g_1(t)},  \,\,\
 t\in[0,1]$
 satisfy  the Hammerstein's system of equations, i.e.
 $$H_kf_{2}=\lambda_1g_{2},\,\,\ H_kg_{2}=\lambda_2f_{2} \eqno(5.5)$$
 where $\lambda_1=\omega(f_1)>0$ , $\lambda_2=\omega(f_2)>0$ and $(f_2,g_2)\in(\mathcal M_0)^{2}$.\\
 On the other hand by Lemma 4.2 there exists  $(f_3,g_3)$ a solution of the
Hammerstain's system of equations:
 $$H_{k}f_{3}=g_{3} , \,\,\,\ H_{k}g_{3}=f_{3}.$$
 But this is contradicts to Proposition 5.2. This completes the
 proof. \\

\textbf{6  Existence of periodic Gibbs Measures for Model (2.1): Case $k=2$.}\\
In this section we construct a function $K(t,u)$ such that
corresponding equation (2.6) has a solution $(f,g)$ with $f\neq
g.$ Put
$$ K_{n}(t,u)=\frac{1-b_{n}c_{n}^{3}\sqrt[n]{u-\frac{1}{2}}\left(\sqrt[n]{(u-\frac{1}{2})^{2}}-4\right)^{2}
\sqrt[n]{t-\frac{1}{2}}}{c_{n}^{2}\left(\sqrt[n]{u-\frac{1}{2}}+2\right)^{2}},
\,\,\,\,\,\ t,u\in[0,1] \eqno(6.1)$$ where
$$ b_{n}=\left(\frac{1}{\sqrt[n]{4}}\right)^{(n-1)}\left(1+\frac{2}{n}\right), \,\,\,\
c_{n}^{3}=\frac{1}{2}\int_{-\frac{1}{2}}^{\frac{1}{2}}\frac{1}{(2+\sqrt[n]{u})^{2}}du.$$

\textbf{Lemma 6.1} For all $t,u \in[0,1],$ the following holds:
$$\lim_{n\rightarrow\infty}K_{n}(t,u)>0.$$

\emph{Proof} It is easy to see
$$\lim_{n\rightarrow\infty}K_{n}(t,u)>0 \,\ \Leftrightarrow \,\ \lim_{n\rightarrow\infty}\left(1-b_{n}c_{n}^{3}
\sqrt[n]{u-\frac{1}{2}}\left(\sqrt[n]{(u-\frac{1}{2})^{2}}-4\right)^{2}
\sqrt[n]{t-\frac{1}{2}}\right)>0.$$
 We have
$$\lim_{n\rightarrow\infty}b_{n}=\lim_{n\rightarrow\infty}\left(\frac{1}{\sqrt[n]{4}}\right)^{(n-1)}\left(1+\frac{2}{n}\right)=\frac{1}{4},$$
$$\lim_{n\rightarrow\infty}c_{n}=\lim_{n\rightarrow\infty}\sqrt[3]{\frac{1}{2}\int_{-\frac{1}{2}}^{\frac{1}{2}}\frac{1}{(2+\sqrt[n]{u})^{2}}du}
\geq\sqrt[3]{\frac{1}{8}}.$$
 Then
$$ \lim_{n\rightarrow\infty}K_{n}(t,u)>0 \,\ \Leftrightarrow \,\ \lim_{n\rightarrow\infty}\left(1-b_{n}c_{n}^{3}\sqrt[n]
{u-\frac{1}{2}}\left(\sqrt[n]{\left(u-\frac{1}{2}\right)^{2}}-4\right)^{2}
\sqrt[n]{t-\frac{1}{2}}\right)\geq$$
$$\geq1-\frac{1}{4}\cdot\frac{1}{8}\cdot9=\frac{23}{32}>0. $$\\

\textbf{Corollary 6.2} \,\ There exists $n_{0}$ such that for
every $n\geq n_{0}$ the function $K_{n_{0}}(t,u)$ is a positive
function.\\

 \emph{Proof} \,\ Straightforward.\\

 \textbf{Theorem 6.3}  The system of
Hammerstain's equation:
$$\int_{0}^{1}K_{n_{0}}(t,u)f^{2}(u)du=g(t), \,\,\,\ \int_{0}^{1}K_{n_{0}}(t,u)g^{2}(u)du=f(t) \eqno(6.2)$$
in the space $(C[0,1])^{2}$ has at least two positive solutions with $f\neq g$.\\

\emph{Proof} Let
$$f^{(n_{0})}_{1}(t)=c_{n_{0}}\left(\sqrt[n_{0}]{t-\frac{1}{2}}+2\right), \,\  g^{(n_{0})}_{1}(t)=1, \,\,\,\ t\in[0,1],$$

Then $(f^{(n_{0})}_{1},g^{(n_{0})}_{1})\in (C[0,1])^{2}$ and positive. \\

(a) Consider the first equation:
$$ \int_{0}^{1}K_{n_{0}}(t,u)f^{2}(u)du=g(t).$$
$$\int_{0}^{1}K_{n_{0}}(t,u)\left(f^{(n_{0})}_{1}(u)\right)^{2}du=
1-\int_{0}^{1}b_{n_{0}}\cdot
c_{n_{0}}^{3}\sqrt[n_{0}]{u-\frac{1}{2}}\left(\sqrt[n_{0}]{\left(u-\frac{1}{2}\right)^{2}-4}\right)^{2}
\sqrt[n_{0}]{t-\frac{1}{2}}du=$$
$$=1-b_{n_{0}}\cdot
c_{n_{0}}^{3}\cdot\sqrt[n_{0}]{t_{1}}\int_{-\frac{1}{2}}^{\frac{1}{2}}\sqrt[n_{0}]{u_{1}}\left(\sqrt[n_{0}]{u_{1}^{2}-4}\right)^{2}du_{1}
=1=g^{(n_{0})}_{1}(t).$$
where $u_{1}=u-\frac{1}{2}$, $t_{1}=t-\frac{1}{2}.$\\

(b) Now we consider the second equation:
 $$\int_{0}^{1}K_{n_{0}}(t,u)\left(g^{(n_{0})}_{1}(u)\right)^{2}du=f^{n_{0}}_{1}(t).$$
$$\int_{0}^{1}K_{n_{0}}(t,u)\left(g^{n_{0}}_{1}(u)\right)^{2}du=\int_{0}^{1}
\frac{1-b_{n_{0}}c_{n_{0}}^{3}\sqrt[n_{0}]{u-\frac{1}{2}}\left(\sqrt[n_{0}]{(u-\frac{1}{2})^{2}}-4\right)^{2}
\sqrt[n_{0}]{t-\frac{1}{2}}}{c_{n_{0}}^{2}\left(\sqrt[n_{0}]{u-\frac{1}{2}}+2\right)^{2}}du=$$
Let $u_{1}=u-\frac{1}{2},$ $t_{1}=t-\frac{1}{2}.$ Then
$$\int_{-\frac{1}{2}}^{\frac{1}{2}}\frac{1}{c^{2}(\sqrt[n_{0}]{u_{1}}+2)^{2}}du_{1}-b_{n_{0}}c_{n_{0}}\sqrt[n_{0}]{t_{1}}\int_{-\frac{1}{2}}^{\frac{1}{2}}
\frac{\sqrt[n_{0}]{u_{1}}\left(\sqrt[n_{0}]{u_{1}}-4\right)^{2}}{(\sqrt[n_{0}]{u_{1}}+2)^{2}}du_{1}=$$
$$\frac{1}{c_{n_{0}}^{2}}\int_{-\frac{1}{2}}^{\frac{1}{2}}\frac{1}{(\sqrt[n_{0}]{u_{1}}+2)^{2}}du_{1}-b_{n_{0}}c_{n_{0}}
\sqrt[n_{0}]{t_{1}}\int_{-\frac{1}{2}}^{\frac{1}{2}}
\sqrt[n_{0}]{u_{1}}\left(\sqrt[n_{0}]{u_{1}}-2\right)^{2}du_{1}=$$
$$=2c_{n_{0}}+4b_{n_{0}}c_{n{0}}\sqrt[n_{0}]{t_{1}}\int_{-\frac{1}{2}}^{\frac{1}{2}}\sqrt[n_{0}]{u_{1}^{2}}du_{1}
=c_{n_{0}}\left(\sqrt[n_{0}]{t-\frac{1}{2}}+2\right)=f^{n_{0}}_{1}(t).$$
 By symmetry of $(f,g)$ we have $\left(g^{n_{0}}_{1}(t),f^{n_{0}}_{1}(t)\right)$ is also
 solution of (6.2).\\
This completes the proof. 

From this we get\\

 \textbf{Theorem 6.4} The model:
$$H(\sigma)=-\frac{1}{\beta}\sum_{<x,y>}ln\left(\frac{1-b_{n_{0}}c_{n_{0}}^{3}\sqrt[n_{0}]{\sigma(x)-\frac{1}{2}}
\left(\sqrt[n_{0}]{(\sigma(x)-\frac{1}{2})^{2}}-4\right)^{2}
\sqrt[n_{0}]{\sigma(y)-\frac{1}{2}}}{c_{n_{0}}^{2}\left(\sqrt[n_{0}]{\sigma(x)-\frac{1}{2}}+2\right)^{2}}\right)$$
 on the Cayley tree $\Gamma^{2}$ has at least two periodic Gibbs
 measures.\\

\textbf{7 Existence of periodic Gibbs Measures for Model (2.1): Case $k=3$.}\\

\textbf{Lemma 7.1} Let $a\in R$. Then for every odd (even)
function $\varphi(x)\in C[0,1]$ the following equation holds:
$$\int_{-a}^{a}\frac{\varphi(x)}{(1+\sin x)^{3}}dx=
-2\int_{0}^{a}\frac{\varphi(x)\sin x(3+\sin^{2}x)}{\cos^{6}x}dx.$$
$$\left(\int_{-a}^{a}\frac{\varphi(x)}{(1+\sin x)^{3}}dx=
2\int_{0}^{a}\frac{\varphi(x)(1+3\sin^{2}x)}{\cos^{6}x}dx
\right).$$\\

\emph{Proof} Let $\varphi(x)$ be odd (the case even is similar)
function
$$\int_{-a}^{a}\frac{\varphi(x)}{(1+\sin x)^{3}}dx=
\int_{0}^{a}\frac{\varphi(x)}{(1+\sin
x)^{3}}dx+\int_{-a}^{0}\frac{\varphi(x)}{(1+\sin x)^{3}}dx=$$
$$\int_{0}^{a}\frac{\varphi(x)}{(1+\sin x)^{3}}dx-
\int_{0}^{a}\frac{\varphi(x)}{(1-\sin x)^{3}}dx=
-2\int_{0}^{a}\frac{\varphi(x)\sin x(3+\sin^{2}x)}{\cos^{6}x}dx. $$\\
Put
$$K(t,u)=\frac{1-\frac{22}{17}\sin\frac{\pi(2t-1)}{3}\sin\frac{\pi(2u-1)}{3}}
{a^{3}(1+\sin\frac{\pi(2u-1)}{3})^{3}}, \,\,\,\,\,\ t,u\in[0,1],
\eqno(7.1)$$ where $a=\sqrt[4]{\frac{198\sqrt{3}}{5\pi}}.$ It is
easy to see that $K(t,u)$ is a positive and continuous function.\\

\textbf{Theorem 7.2}  The system of Hammerstain's equations
$$\int_{0}^{1}K(t,u)f^{3}(u)du=g(t), \,\,\,\ \int_{0}^{1}K(t,u)g^{3}(u)du=f(t), \eqno(7.2)$$
in the space $(C[0,1])^{2}$ has at least two positive solutions with $f\neq g$.\\

\emph{Proof} (a) Denote
$$f_{1}(t)=a\left(1+\sin\frac{\pi(2t-1)}{3}\right), g_{1}(t)=1,  \,\,\
t\in[0,1],$$ where $a=\sqrt[4]{\frac{198\sqrt{3}}{5\pi}}.$ Then
$(f_{1},g_{1})\in (C[0,1])^{2}$ and the functions $f_{1}$ and
$g_{1}$ are positive.
 Consider the first equation of (7.2)
$$\int_{0}^{1}K(t,u)f_{1}^{3}(u)du=1-\frac{22}{17}\sin\frac{\pi(2t-1)}{3}\int_{0}^{1}\sin\frac{\pi(2u-1)}{3}du=1.$$
(b) Now we check the second equation.
$$\int_{0}^{1}K(t,u)g^{3}_{1}(u)du=\int_{0}^{1}\frac{1-\frac{22}{17}\sin\frac{\pi(2t-1)}{3}\sin\frac{\pi(2u-1)}{3}}
{a^{3}(1+\sin\frac{\pi(2u-1)}{3})^{3}}du.$$
 Let $t_{1}=\frac{\pi}{3}(2t-1)$, $u_{1}=\frac{\pi}{3}(2u-1).$ Then
 $$\int_{-\frac{\pi}{2}}^{\frac{\pi}{2}}K(t_{1},u_{1})g_{1}^{3}(u_{1})du_{1}=\frac{3}{2a^{3}\pi}\left(\int_{-\frac{\pi}{3}}^{\frac{\pi}{3}}
 \frac{1-\frac{22}{17}\sin t_{1}\sin u_{1}}{(1+\sin u_{1})^{3}}\right)du_{1}=$$
$$=\frac{3}{2a^{3}\pi}\left(\int_{-\frac{\pi}{3}}^{\frac{\pi}{3}}
 \frac{1}{(1+\sin u_{1})^{3}}du_{1}-
 \frac{22}{17}\sin t_{1}\int_{-\frac{\pi}{3}}^{\frac{\pi}{3}}
 \frac{\sin u_{1}}{(1+\sin u_{1})^{3}}\right)du_{1}.$$
 By Lemma 7.1 LHS of this equality is
  $$\frac{3}{2a^{3}\pi}\left[\int_{0}^{\frac{\pi}{3}}\frac{1+3\sin^{2}u_{1}}{\cos^{6}u_{1}}du_{1}
  +\frac{22}{17}\sin t_{1}\int_{0}^{\frac{\pi}{3}}\frac{\sin^{2}u_{1}(3+\sin^{2}u_{1})}{\cos^{6}u_{1}}du_{1}\right]=$$
  $$=\frac{3}{a^{3}\pi}\left[\int_{0}^{\sqrt{3}}(1+4y^{2})(1+y^{2})dy+
  \frac{22}{17}\sin t_{1}\int_{0}^{\sqrt{3}}y^{2}(3+4y^{2})dy\right]=$$
  $$=\frac{198\sqrt{3}}{a^{3}\pi}(1+\sin t_{1})=\frac{198\sqrt{3}}{a^{3}\pi}\left(1+\sin\frac{\pi(2t-1)}{3}\right)=f(t).$$
By symmetry of $(f_{1},g_{1})$ we have $(g_{1}(t),f_{1}(t))$ is also solution to (6.2).\\
  This completes the proof. \\

\textbf{Theorem 7.3} The model:
$$H(\sigma)=-\frac{1}{\beta}\sum_{<x,y>}ln\left(\frac{1-\frac{22}{17}\sin\frac{\pi(2\sigma(x)-1)}{3}\sin\frac{\pi(2\sigma(y)-1)}{3}}
{a^{3}(1+\sin\frac{\pi(2\sigma(x)-1)}{3})^{3}}\right)$$
 on the Cayley tree $\Gamma^{3}$ has at least two periodic Gibbs
 measures.\\

\textbf{8  Existence of periodic Gibbs Measures for Model (2.1): Case $k\geq4$.}\\

 Denote
$$c_{k}=\frac{2\left(1-(\frac{1}{3})^{k-1}\right)}
{\frac{k-1}{k-2}\left(1-(\frac{1}{3})^{k-2}\right)-2\left(1-(\frac{1}{3})^{k-1}\right)}.
\eqno(8.1)$$\\

 \textbf{Lemma 8.1} For every $k\in N ,\,\  k\geq4$  the following
 inequality holds: $|c_{k}|<4$.\\

 \emph{Proof} \,\ For $k\geq4$ we have
$$|c_{k}|=\frac{2\left(1-(\frac{1}{3})^{k-1}\right)}
{\frac{k-3}{k-2}+(\frac{1}{3})^{k-1}(\frac{k+1}{k-2})}<\frac{2}
{\frac{k-3}{k-2}+(\frac{1}{3})^{k-1}(\frac{k+1}{k-2})}<\frac{2(k-2)}{k-3},$$
Thus
$$|c_{k}|<\frac{2(k-2)}{k-3}=2+\frac{2}{k-3}\leq4.  \eqno(8.2)$$
Hence $|c_{k}|<4$ for $k\geq4$. \\

For each $k\geq4$ , $a>0$ we define the continuous function
$$K(t,u,k)=\frac{1+c_{k}(t-\frac{1}{2})(u-\frac{1}{2})}{a^{k}(u+\frac{1}{2})^{k}},
 \,\,\,\,\,\,\,\ t,u\in[0,1].$$
By the inequality (8.2) it follows that the function
$K(t,u,k)$ is positive.\\

\textbf{Theorem 8.2} For each $k\geq4$ the Hammerstein's system of
equations:
$$\int_{0}^{1}K(t,u,k)f^{k}(u)du=g(t) ,\,\ \int_{0}^{1}K(t,u,k)g^{k}(u)du=f(t) \eqno(8.3)$$
in $(C[0,1])^{2}$ have at least two positive solutions with $f\neq g$.\\

 \emph{Proof}  Let $k\geq4$. Define the positive
 continuous functions $f_{1}(t)$, $g_{1}(t)$ on $[0,1]$ by the
 equality
 $$f_{1}(t)=a\left(t+\frac{1}{2}\right), \,\,\,\ g_{1}(t)=1$$
 where
 $$a=a(k)=\sqrt[k+1]{\frac{2^{k-1}}{k-1}\left(1-\left(\frac{1}{3}
 \right)^{k-1}\right)} , \,\,\,\,\,\ k\geq4.$$
 It is easy to see that $a>0$. We shall show that $(f_{1},g_{1})$ is a solution to the Hammerstein's
 system of equations (8.2).\\

 We shall check the first equation.
 $$\int_{0}^{1}K(t,u,k)f_{1}^{k}(u)du=\int_{0}^{1}\frac{1+c_{k}\left(t-\frac{1}{2}\right)
\left(u-\frac{1}{2}\right)}{a^{k}\left(u+\frac{1}{2}\right)^{k}}
\left(a\left(u+\frac{1}{2}\right)\right)^{k}du=$$
$$=\int_{0}^{1}\left(1+c_{k}\left(t-\frac{1}{2}\right)\left(u-\frac{1}{2}\right)\right)du=1+c_{k}
t_{1}\int_{-\frac{1}{2}}^{\frac{1}{2}}u_{1}du_{1}=1.$$ Where
$t_{1}=t-\frac{1}{2}$ and $u_{1}=u-\frac{1}{2}$. Hence
$$\int_{0}^{1}K(t,u,k)f_{1}^{k}(u)du=g_{1}(t).$$
Now we shall check the second equation.
$$\int_{0}^{1}K(t,u,k)g_{1}^{k}(u)du=\int_{0}^{1}\left(\frac{1+c_{k}
(t-\frac{1}{2})(u-\frac{1}{2})}{a^{k}(u+\frac{1}{2})^{k}}\right)du=$$
$$\int_{0}^{1}\frac{1}{a^{k}(u+\frac{1}{2})^{k}}du+\int_{0}^{1}
\frac{c_{k}(t-\frac{1}{2})(u-\frac{1}{2})}{a^{k}(u+\frac{1}{2})^{k}}du=
\frac{1}{a^{k}}\left(\int_{-\frac{1}{2}}^{\frac{1}{2}}\frac{1}
{(u_{1}+1)^{k}}du_{1}+c_{k}t_{1}
\int_{-\frac{1}{2}}^{\frac{1}{2}}\frac{u_{1}}{(u_{1}+1)^{k}}du_{1}\right)=$$
 (where $t_{1}=t-\frac{1}{2}$, $u_{1}=u-\frac{1}{2}$)
$$=\frac{1}{a^{k}}\frac{2^{k-1}}{k-1}\left(1-\left(\frac{1}{3}\right)^{k-1}\right)
+\frac{c_{k}t_{1}}{a^{k}}\left[\frac{2^{k-2}}{k-2}\left(1-\left(\frac{1}{3}\right)^{k-2}\right)
-\frac{2^{k-1}}{k-1}\left(1-\left(\frac{1}{3}\right)^{k-1}\right)\right]=$$
$$=a+at_{1}\frac{2\left(1-(\frac{1}{3})^{k-1}\right)}
{\frac{k-1}{k-2}\left(1-(\frac{1}{3})^{k-2}\right)-
2\left(1-(\frac{1}{3})^{k-1}\right)}\left(\frac{(k-1)
\left(1-(\frac{1}{3})^{k-2}\right)-2(k-2)\left(1-(\frac{1}{3})^{k-1}
\right)}{2(k-2)\left(1-(\frac{1}{3})^{k-1}\right)}\right)=$$
$$=a+at_{1}=a\left(t+\frac{1}{2}\right)=f_{1}(t).$$
 Moreover, $\left(g_{1}(t),f_{1}(t)\right)$ is also solution to (8.3). \\

\textbf{Theorem 8.3} Let $k\geq4.$ The model
$$H(\sigma)=-\frac{1}{\beta}\sum_{<x,y>}ln\left(\frac{1+c_{k}(\sigma(x)-\frac{1}{2})(\sigma(y)-\frac{1}{2})}{a^{k}
(\sigma(x)+\frac{1}{2})^{k}}\right)$$
 on the Cayley tree $\Gamma^{k}$ has at least two periodic Gibbs
 measures.\\

 \textbf{9  Existence of four periodic Gibbs Measures for Model (2.1).}\\

Denote
$$c_{ij}(m)=\frac{1}{m+2(i-1)+2(j-1)}, \,\,\ (n,m,p)\in
{N\times N\times N_{0}}, \,\,\  1\leq i,j\leq n,$$
$A_{n}^{(m,p)}=\left(\frac{c_{ij}(m)}{4^{p+j+i-2}}\right)_{n}$ be $n\times n$ square matrix.\\
If $n\in\{2,3\}$ then it's easy to check  $\det
(A_{n}^{(m,p)})\neq0.$\\

Put
    $$\begin{pmatrix}
     a_{11}\\
     a_{12}\\
     a_{13}
     \end{pmatrix}=\left(A_{3}^{(1,0)}\right)^{-1}
     \begin{pmatrix}0\\
     \frac{1}{6}\\
     0\end{pmatrix}, \,\,\ \begin{pmatrix}
     a_{21}\\
     a_{22}\\
     a_{23}
     \end{pmatrix}=\left(A_{3}^{(3,1)}\right)^{-1}
     \begin{pmatrix}0\\
     0\\
     \frac{1}{20}\end{pmatrix}$$
and
    $$\begin{pmatrix}
    b_{11}\\
    b_{12}
    \end{pmatrix}=\left(A_{2}^{(5,2)}\right)^{-1}
    \begin{pmatrix}1\\
    0\end{pmatrix}, \,\,\
    \begin{pmatrix}b_{21}\\
    b_{22}
    \end{pmatrix}=\left(A_{2}^{(7,3)}\right)^{-1}
    \begin{pmatrix}0\\
    1
    \end{pmatrix},$$
where $\left(A_{n}^{(m,p)}\right)^{-1}$ is inverse of
    $A_{n}^{(m,p)}.$\\
So we define following functions:
    $$\psi_{1}(u)=a_{11}+a_{12}u^{2}+a_{13}u^{4}, \,\,\ \psi_{2}(u)=a_{21}u^{2}+a_{22}u^{4}+a_{23}u^{6},$$
    $$\psi_{3}(u)=b_{11}u+b_{12}u^{3}, \,\,\
    \psi_{4}(u)=b_{21}u^{3}+b_{22}u^{5}.$$
Finally
$$K_{1}(t,u;k)=\psi_{1}(u)\left(\sqrt[k]{20t^{4}+\frac{3}{4}}-1\right)+\psi_{2}(u)\left(\sqrt[k]{6t^{2}+\frac{1}{2}}-1\right),$$
$$K_{2}(t,u;k)=\psi_{3}(u)\left(\sqrt[k]{t^{3}+1}-1\right)+\psi_{4}(u)\left(\sqrt[k]{t^{5}+1}-1\right),$$
$$\widetilde{K}(t,u;k)=1+K_{1}(t,u;k)+K_{2}(t,u;k).$$

\textbf{Remark 9.1} \,\ There exist $k_{0}\in N$ such that for all
$k\geq k_{0}$ the following inequality holds
 $$\widetilde{K}\left(t-\frac{1}{2},u-\frac{1}{2};k\right)>0,
\,\,\,\ (t,u)\in [0,1]^{2}.$$

\emph{Proof} \,\ It is sufficient to show:
$$\lim_{k\rightarrow\infty}\widetilde{K}\left(t-\frac{1}{2},u-\frac{1}{2};k\right)>0, \,\,\,\ (t,u)\in [0,1]^{2}.$$
Let $\gamma:[0,1]\rightarrow[m,M]$ be a function, $m>0.$\\
 We have
$$0=\lim_{k\rightarrow\infty}(\sqrt[k]{m}-1)\leq\lim_{k\rightarrow\infty}(\sqrt[k]{\gamma(t)}-1)
\leq\lim_{k\rightarrow\infty}(\sqrt[k]{M}-1)=0.$$
 Hence
$$\lim_{k\rightarrow\infty}(\sqrt[k]{\gamma(t)}-1)=0 \,\ \Rightarrow \lim_{k\rightarrow\infty}K_{i}
\left(t-\frac{1}{2},u-\frac{1}{2};k\right)=0, \,\,\,\
i\in\{1,2\}$$
 and
 $$\lim_{k\rightarrow\infty}\widetilde{K}\left(t-\frac{1}{2},u-\frac{1}{2};k\right)=1>0.$$
This completes the proof. \\

\textbf{Lemma 9.2} \,\ If $k\in\{1,2\}$, then
$$(i)\,\
\int_{-\frac{1}{2}}^{\frac{1}{2}}\psi_{1}(u)u^{2k}du=\frac{1}{12}\left(1+(-1)^{k+1}\right),
\,\,\,\,\,\,\,\,\,\,\,\,\,\,\,\,\,\,\,\,\,\,\,\,\,\,\,\,\,\,\,\,\,\,\,\,\,\,\,\,\,\,\,\,\,\,\,\,\,\,\,\,\,\,\,\,\,\,\,\,\,\,\,\,\,\,\,\,\,\,\,\,\,\,\,\,\,\,\,\,\,\,\,\,\,\,\,\,\,\,\,\,\,
\,\,\,\,\,\,\,\,\,\,\,\,\,\,\,\,\,\  $$
$$(ii)\,\
\int_{-\frac{1}{2}}^{\frac{1}{2}}\psi_{2}(u)u^{2k}du=\frac{1}{12}\left(1+(-1)^{k}\right),
\,\,\,\,\,\,\,\,\,\,\,\,\,\,\,\,\,\,\,\,\,\,\,\,\,\,\,\,\,\,\,\,\,\,\,\,\,\,\,\,\,\,\,\,\,\,\,\,\,\,\,\,\,\,\,\,\,\,\,\,\,\,\,\,\,\,\,\,\,\,\,\,\,\,\,\,\,\,\,\,\,\,\,\,\,\,\,\,\,\,\,\,\,\,\,
\,\,\,\,\,\,\,\,\,\,\,\,\,\,\,\,\,\,\ $$
$$(iii)\,\ \int_{-\frac{1}{2}}^{\frac{1}{2}}\psi_{3}(u)u^{2k+1}du=\frac{1}{12}\left(1+(-1)^{k}\right), \,\,\,\,\,\,\,\,\,\,\,\,\,\,\,\,\,\,\,\,\,\,\,\,\,\,\,\,\,\,\,\,\,\,\,\,\,\,\,\,\,\,\,\,\,\,\,\,\,\,\,\,\,\,\,\,\,\,\,\,\,\,\,\,\,\,\,\,\,\,\,\,\,\,\,\,\,\,\,\,\,\,\,\,\,\,\,\,\,\,\,\,\,\,\,
\,\,\,\,\,\,\,\,\,\,\,\,\,\,\,\,\ $$
$$(iv)\,\
\int_{-\frac{1}{2}}^{\frac{1}{2}}\psi_{4}(u)u^{2k+1}du=\frac{1}{12}\left(1+(-1)^{k+1}\right).
\,\,\,\,\,\,\,\,\,\,\,\,\,\,\,\,\,\,\,\,\,\,\,\,\,\,\,\,\,\,\,\,\,\,\,\,\,\,\,\,\,\,\,\,\,\,\,\,\,\,\,\,\,\,\,\,\,\,\,\,\,\,\,\,\,\,\,\,\,\,\,\,\,\,\,\,\,\,\,\,\,\,\,\,\,\,\,\,\,\,\,\,\,\,\,
\,\,\,\,\,\,\,\,\,\,\,\ $$

\emph{Proof} \,\ For $k\in \{1,2\}$
$$(i)\,\,\,\ \int_{-\frac{1}{2}}^{\frac{1}{2}}\psi_{1}(u)u^{2k}du=
a_{11}\int_{-\frac{1}{2}}^{\frac{1}{2}}u^{2k}du+a_{12}\int_{-\frac{1}{2}}^{\frac{1}{2}}u^{2(k+1)}du
+a_{13}\int_{-\frac{1}{2}}^{\frac{1}{2}}u^{2(k+2)}du=
\,\,\,\,\,\,\,\,\,\,\,\,\ $$
$$=\frac{a_{11}}{4^{k}(2k+1)}+\frac{a_{12}}{4^{k+1}(2k+3)}+\frac{a_{13}}{4^{k+2}(2k+5)}=$$
$$=a_{11}\times c_{k+1,1}^{(0)}(1)+a_{12}\times
c_{k+1,2}^{(0)}(1)+a_{13}\times
c_{k+1,3}^{(0)}(1)=\frac{1}{12}\left(1+(-1)^{k+1}\right).$$
$$(ii) \,\ \int_{-\frac{1}{2}}^{\frac{1}{2}}\psi_{2}(u)u^{2k}du=
a_{21}\int_{-\frac{1}{2}}^{\frac{1}{2}}u^{2(k+1)}du+a_{22}
\int_{-\frac{1}{2}}^{\frac{1}{2}}u^{2(k+2)}du+a_{23}\int_{-\frac{1}{2}}^{\frac{1}{2}}u^{2(k+3)}du=\,\,\,\
$$
$$=\frac{a_{21}}{4^{k}(2k+1)}+\frac{a_{22}}{4^{k+1}(2k+3)}+\frac{a_{23}}{4^{k+2}(2k+5)}=$$
$$=a_{21}\times c_{k+1,1}^{(1)}(3)+a_{22}\times
c_{k+1,2}^{(1)}(3)+a_{23}\times
c_{k+1,3}^{(1)}(3)=\frac{1}{12}\left(1+(-1)^{k}\right).$$
$$(iii) \,\,\ \int_{-\frac{1}{2}}^{\frac{1}{2}}\psi_{3}(u)u^{2k+1}du=
b_{11}\int_{-\frac{1}{2}}^{\frac{1}{2}}u^{2(k+1)}du+b_{12}\int_{-\frac{1}{2}}^{\frac{1}{2}}u^{2(k+2)}du=
\,\,\,\,\,\,\,\,\,\,\,\,\,\,\,\,\,\,\,\,\,\,\,\,\,\,\,\,\,\,\,\,\,\,\,\,\,\,\,\,\,\,\,\,\,\,\,\,\,\,\,\,\
$$
$$=b_{11}\times c_{k,1}^{(2)}(5)+b_{12}\times
c_{k,2}^{(2)}(5)=\frac{1}{12}\left(1+(-1)^{k}\right)=\frac{1}{12}\left(1+(-1)^{k}\right).$$
$$(iv) \,\,\ \int_{-\frac{1}{2}}^{\frac{1}{2}}\psi_{4}(u)u^{2k+1}du=
b_{21}\int_{-\frac{1}{2}}^{\frac{1}{2}}u^{2(k+2)}du+b_{22}\int_{-\frac{1}{2}}^{\frac{1}{2}}u^{2(k+3)}du=
\,\,\,\,\,\,\,\,\,\,\,\,\,\,\,\,\,\,\,\,\,\,\,\,\,\,\,\,\,\,\,\,\,\,\,\,\,\,\,\,\,\,\,\,\,\,\,\,\,\,\,\
$$
$$=b_{21}\times c_{k,1}^{(3)}(7)+b_{22}\times
c_{k,2}^{(3)}(7)=\frac{1}{12}\left(1+(-1)^{k+1}\right).  $$

\textbf{Lemma 9.3} \,\ The function $\varphi_{0}(u)=1$ is a fixed
point of the operator $H_{k}:$
$$(H_{k}f)(t)=\int_{0}^{1}\widetilde{K}\left(t-\frac{1}{2},u-\frac{1}{2};k\right)f^{k}(u)du, \,\,\ k\geq 2
\eqno(9.1)$$

 \emph{Proof} \,\ Let $u_{1}=u-\frac{1}{2}, \,\
 v_{1}=v-\frac{1}{2}$  then
 $$(H_{k}\varphi_{0})\left(t-\frac{1}{2}\right)=\int_{0}^{1}\widetilde{K}\left(t-\frac{1}{2},u-\frac{1}{2};k\right)du=
 \int_{-\frac{1}{2}}^{\frac{1}{2}}\widetilde{K}(t_{1},u_{1};k)du_{1}=$$
 $$=\int_{-\frac{1}{2}}^{\frac{1}{2}}\left[1+K_{1}(t_{1},u_{1};k)+K_{2}(t_{1},u_{1};k)\right]du_{1}=
 1+\int_{-\frac{1}{2}}^{\frac{1}{2}}K_{1}(t_{1},u_{1};k)du_{1}+\int_{-\frac{1}{2}}^{\frac{1}{2}}K_{2}
 (t_{1},u_{1};k)du_{1}.$$
 Now we'll prove the following
 $$\int_{-\frac{1}{2}}^{\frac{1}{2}}K_{i}(t_{1},u_{1};k)du_{1}=0, \,\,\,\ i\in\{1,2\}. \eqno(9.2)$$
Case: $i=1$
 $$\int_{-\frac{1}{2}}^{\frac{1}{2}}K_{1}(t_{1},u_{1};k)du_{1}=\int_{-\frac{1}{2}}^{\frac{1}{2}}
 \left[\psi_{1}(u_{1})\left(\sqrt[k]{20t_{1}^{4}+\frac{3}{4}}-1\right)du_{1}+
 \psi_{2}(u_{1})\left(\sqrt[k]{6t_{1}^{2}+\frac{1}{2}}-1\right)\right]du_{1}=$$
$$=\left(\sqrt[k]{20t_{1}^{4}+\frac{3}{4}}-1\right)\left(a_{11}+a_{12}
\int_{-\frac{1}{2}}^{\frac{1}{2}}u^{2}du+a_{13}\int_{-\frac{1}{2}}^{\frac{1}{2}}u^{4}du\right)+$$
$$+\left(\sqrt[k]{6t_{1}^{2}+\frac{1}{2}}-1\right)\left(a_{21}\int_{-\frac{1}{2}}^{\frac{1}{2}}u^{2}du+a_{22}
\int_{-\frac{1}{2}}^{\frac{1}{2}}u^{4}du+a_{23}\int_{-\frac{1}{2}}^{\frac{1}{2}}u^{6}du\right)=$$
$$=\left(\sqrt[k]{20t_{1}^{4}+\frac{3}{4}}-1\right)\left(a_{11}+\frac{a_{12}}{3\cdot4}+\frac{a_{13}}{5\cdot4^{2}}\right)+
\left(\sqrt[k]{6t_{1}^{2}+\frac{1}{2}}-1\right)\left(\frac{a_{21}}{3\cdot4}+\frac{a_{22}}{5\cdot4^{2}}+\frac{a_{23}}{7\cdot4^{k+2}}\right)=$$
$$=\left(\sqrt[k]{20t_{1}^{4}+\frac{3}{4}}-1\right)\left(a_{11}\times c_{1,1}^{(0)}(1)+a_{12}\times
c_{1,2}^{(0)}(1)+a_{13}\times c_{1,3}^{(0)}(1)\right)+$$
$$+\left(\sqrt[k]{6t_{1}^{2}+\frac{1}{2}}-1\right)\left(a_{21}\times c_{2,1}^{(1)}(3)+a_{22}\times
c_{2,2}^{(1)}(3)+a_{23}\times c_{2,3}^{(1)}(3)\right)=0.$$
Case:\,\ $i=2$
 $$\int_{-\frac{1}{2}}^{\frac{1}{2}}K_{2}(t_{1},u_{1};k)=\int_{-\frac{1}{2}}^{\frac{1}{2}}\psi_{3}(u_{1})
 \left(\sqrt[k]{t^{3}+1}\right)du_{1}+\int_{-\frac{1}{2}}^{\frac{1}{2}}\psi_{4}(u_{1})
 \left(\sqrt[k]{t^{5}+1}\right)du_{1}$$
It's easy to check for $j\in\{3,4\}$ the functions
$\psi_{j}(u_{1})$ is odd, i.e:
$$\int_{-\frac{1}{2}}^{\frac{1}{2}}\psi_{j}(u_{1})du_{1}=0 \,\ \Rightarrow \,\ \int_{-\frac{1}{2}}^{\frac{1}{2}}K_{2}(t_{1},u_{1};k)du_{1}=0.$$
Thus we have proved
$$\left(H_{k}\varphi_{0}\right)(t)=1+\int_{-\frac{1}{2}}^{\frac{1}{2}}K_{1}(t_{1},u_{1};k)du_{1}+
\int_{-\frac{1}{2}}^{\frac{1}{2}}K_{2}(t_{1},u_{1};k)du_{1}=1.$$
This completes the proof. 

Denote
$$f_{1}(u)=\sqrt[k]{6u^{2}+\frac{1}{2}}, \,\,\
f_{2}(u)=\sqrt[k]{20u^{4}+\frac{3}{4}}, \,\,\
g_{1}(u)=\sqrt[k]{u^{3}+1}, \,\,\ g_{2}(u)=\sqrt[k]{u^{5}+1}.$$

 \textbf{Theorem 9.4} \,\ For all $k\geq k_{0}$ the Hammerstein's
 system of equations:
 $$\int_{0}^{1}\widetilde{K}\left(t-\frac{1}{2},u-\frac{1}{2};k\right)f^{k}(u)du=g(t), \,\,\
 \int_{0}^{1}\widetilde{K}\left(t-\frac{1}{2},u-\frac{1}{2};k\right)g^{k}(u)du=f(t) \eqno(9.3)$$
in $\left(C[0;1]\right)^{2}$ have at least four positive
solutions with $f\neq g.$\\

\emph{Proof} \,\ We'll show
$$\left(f_{1}\left(u-\frac{1}{2}\right),\,\ f_{2}\left(u-\frac{1}{2}\right)\right), \,\ \left(f_{2}\left(u-\frac{1}{2}\right),\,\
f_{1}\left(u-\frac{1}{2}\right)\right)$$
and
$$\left(g_{1}\left(u-\frac{1}{2}\right),\,\ g_{2}\left(u-\frac{1}{2}\right)\right), \,\
\left(g_{2}\left(u-\frac{1}{2}\right),\,\
g_{1}\left(u-\frac{1}{2}\right)\right)$$
 are solutions to the system of equations (9.3).\\

  At first we'll prove $\left(f_{1}(u-\frac{1}{2}),\,\
f_{2}(u-\frac{1}{2})\right)$ is a solution to equation (9.3).\\
 Let
  $u-\frac{1}{2}=u_{1},$\,\ $t-\frac{1}{2}=t_{1}.$ \,\ Then
  $$\left(H_{k}f_{i}\right)(t)=\int_{0}^{1}\widetilde{K}\left(t-\frac{1}{2},u-\frac{1}{2};k\right)f_{i}^{k}\left(u-\frac{1}{2}\right)du=
  \int_{0}^{1}\widetilde{K}(t_{1},u_{1};k)f_{i}^{k}(u_{1})du_{1}=$$
$$=\int_{0}^{1}\left[1+K_{1}(t_{1},u_{1};k)+K_{2}(t_{1},u_{1};k)\right]f_{i}^{k}(u_{1})du_{1}.$$
It's easy to see that
 $$K_{2}(t_{1},-u_{1};k)=-K_{2}(t_{1},u_{1};k), \,\,\,\ f_{i}(u_{1})=f_{i}(-u_{1}), \,\,\,\ i\in\{1,2\}.$$
Hence
 $$K_{2}(t_{1},-u_{1};k)f_{i}(u_{1})=-K_{2}(t_{1},u_{1};k)f_{i}(u_{1}) \,\ \Rightarrow \,\
 \int_{-\frac{1}{2}}^{\frac{1}{2}}K_{2}(t_{1},u_{1};k)f_{i}(u_{1})du_{1}=0. $$
Thus
$$\left(H_{k}f_{i}\right)\left(t-\frac{1}{2}\right)=\left(H_{k}f_{i}\right)(t_{1})=\int_{-\frac{1}{2}}^{\frac{1}{2}}
\left[1+K_{1}(t_{1},u_{1};k)\right]f_{i}^{k}(u_{1})du_{1}.
\eqno(9.4)$$ Case: $i=1$
$$\left(H_{k}f_{1}\right)\left(t-\frac{1}{2}\right)=\left(H_{k}f_{1}\right)(t_{1})=\int_{-\frac{1}{2}}^{\frac{1}{2}}
\left[1+K_{1}(t_{1},u_{1};k)\right]\left(\sqrt[k]{6u_{1}^{2}+\frac{1}{2}}\right)^{k}du_{1}=$$
$$=\int_{-\frac{1}{2}}^{\frac{1}{2}}\left(6u_{1}^{2}+\frac{1}{2}\right)du_{1}+
\int_{-\frac{1}{2}}^{\frac{1}{2}}K_{1}(t_{1},u_{1};k)\left(6u_{1}^{2}+\frac{1}{2}\right)du_{1}=$$
$$1+6\int_{-\frac{1}{2}}^{\frac{1}{2}}K_{1}(t_{1},u_{1};k)u^{2}_{1}du_{1}+
\frac{1}{2}\int_{-\frac{1}{2}}^{\frac{1}{2}}K_{1}(t_{1},u_{1};k)du_{1}=$$
By (9.2) we get
$$=1+6\int_{-\frac{1}{2}}^{\frac{1}{2}}\left[\psi_{1}(u_{1})\left(\sqrt[k]{20t_{1}^{4}+\frac{3}{4}}-1\right)+\psi_{2}(u_{1})
\left(\sqrt[k]{6t_{1}^{2}+\frac{1}{2}}-1\right)\right]u_{1}^{2}du_{1}=$$
$$=1+6\left(\sqrt[k]{20t_{1}^{4}+\frac{3}{4}}-1\right)
\int_{-\frac{1}{2}}^{\frac{1}{2}}\psi_{1}(u_{1})u_{1}^{2}du_{1}+6\left(\sqrt[k]{6t_{1}^{2}+\frac{1}{2}}-1\right)
\int_{-\frac{1}{2}}^{\frac{1}{2}}\psi_{2}(u_{1})u_{1}^{2}du_{1}.$$
By Lemma 9.2
$$\int_{-\frac{1}{2}}^{\frac{1}{2}}\psi_{1}(u_{1})u_{1}^{2}du_{1}=\frac{1}{6},\,\,\
\int_{-\frac{1}{2}}^{\frac{1}{2}}\psi_{2}(u_{1})u_{1}^{2}du_{1}=0.
\eqno(9.5)$$ By (9.2) and (9.5) we obtain
$$(H_{k}f_{1})\left(t-\frac{1}{2}\right)=1+\left(\sqrt[k]{20t_{1}^{4}+\frac{3}{4}}-1\right)=
\sqrt[k]{20t_{1}^{4}+\frac{3}{4}}=f_{2}(t_{1})=f_{2}\left(t-\frac{1}{2}\right).$$
Case: $i=2$
$$\left(H_{k}f_{2}\right)\left(t-\frac{1}{2}\right)=\left(H_{k}f_{2}\right)(t_{1})=\int_{-\frac{1}{2}}^{\frac{1}{2}}
\left[1+K_{1}(t_{1},u_{1};k)\right]f_{2}^{k}(u_{1})du_{1}=$$
$$=\int_{-\frac{1}{2}}^{\frac{1}{2}}\left[1+K_{1}(t_{1},u_{1};k)\right]\left(20u^{4}_{1}+\frac{1}{2}\right)du_{1}=$$
$$=\frac{1}{2}\int_{-\frac{1}{2}}^{\frac{1}{2}}K_{1}(t_{1},u_{1};k)du_{1}+20\int_{-\frac{1}{2}}^{\frac{1}{2}}K_{1}(t_{1},u_{1};k)u^{2}_{1}du_{1}+
\int_{-\frac{1}{2}}^{\frac{1}{2}}\left(20u^{4}_{1}+\frac{1}{2}\right)du_{1}=$$
By  (9.2)
$$=1+20\int_{-\frac{1}{2}}^{\frac{1}{2}}K_{1}(t_{1},u_{1};k)u^{4}_{1}du_{1}=$$
$$=1+20\left(\sqrt[k]{20t_{1}^{4}+\frac{3}{4}}-1\right)\int_{-\frac{1}{2}}^{\frac{1}{2}}\psi_{1}(u_{1})u_{1}^{4}du_{1}
+20\left(\sqrt[k]{6t_{1}^{2}+\frac{1}{2}}-1\right)
\int_{-\frac{1}{2}}^{\frac{1}{2}}\psi_{2}(u_{1})u_{1}^{4}du_{1}.$$
By Lemma 9.2
$$\int_{-\frac{1}{2}}^{\frac{1}{2}}\psi_{1}(u_{1})u_{1}^{4}du_{1}=0, \,\,\
\int_{-\frac{1}{2}}^{\frac{1}{2}}\psi_{2}(u_{1})u_{1}^{4}du_{1}=\frac{1}{20}.$$
Then
$$\left(H_{k}f_{2}\right)\left(t-\frac{1}{2}\right)=1+\left(\sqrt[k]{6t_{1}^{2}+\frac{1}{2}}-1\right)=
\sqrt[k]{6t_{1}^{2}+\frac{1}{2}}=f_{1}(t_{1})=f_{1}\left(t-\frac{1}{2}\right).$$
 By symmetry of $(f_{1},f_{2})$ we have $(f_{2},f_{1})$
is also solution to equation (9.3).\\

Now we'll prove
$\left(g_{1}(u-\frac{1}{2}),g_{2}(u-\frac{1}{2})\right)$ is a
solution to equation (9.3).\\
For $i\in\{1,2\}$
$$(H_{k}g_{i})\left(t-\frac{1}{2}\right)=\int_{0}^{1}\widetilde{K}\left(t-\frac{1}{2},u-\frac{1}{2};k\right)
g_{i}^{k}\left(u-\frac{1}{2}\right)du
=\int_{-\frac{1}{2}}^{\frac{1}{2}}\widetilde{K}(t_{1},u_{1};k)(1+u_{1}^{2i+1})du_{1},$$
where $u_{1}=u-\frac{1}{2}, \,\ t_{1}=t-\frac{1}{2}.$\,\
Then
$$(H_{k}g_{i})\left(t-\frac{1}{2}\right)=\int_{-\frac{1}{2}}^{\frac{1}{2}}\widetilde{K}(t_{1},u_{1};k)du_{1}+
\int_{-\frac{1}{2}}^{\frac{1}{2}}\widetilde{K}(t_{1},u_{1};k)u_{1}^{2i+1}du_{1}=$$
By Lemma 9.3
$$=(H_{k}g_{i})\left(t-\frac{1}{2}\right)=1+\int_{-\frac{1}{2}}^{\frac{1}{2}}\widetilde{K}(t_{1},u_{1};k)u_{1}^{2i+1}du_{1}=$$
$$=1+\int_{-\frac{1}{2}}^{\frac{1}{2}}[1+K_{1}(t_{1},u_{1};k)+K_{2}(t_{1},u_{1};k)]u_{1}^{2i+1}du_{1}.$$
Hence
$$(H_{k}g_{i})\left(t-\frac{1}{2}\right)=1+\int_{-\frac{1}{2}}^{\frac{1}{2}}[1+K_{1}(t_{1},u_{1};k)]u_{1}^{2i+1}du_{1}
+\int_{-\frac{1}{2}}^{\frac{1}{2}}K_{2}(t_{1},u_{1};k)u_{1}^{2i+1}du_{1}
\eqno(9.6)$$ One can easily check that
$$K_{1}(t_{1},-u_{1};k)=K_{1}(t_{1},u_{1};k)\,\,\ \Rightarrow \,\ K_{1}(t_{1},-u_{1};k)(-u_{1}^{2i+1})=-K_{1}(t_{1},u_{1};k)u_{1}^{2i+1},$$
then
$$\int_{-\frac{1}{2}}^{\frac{1}{2}}K_{1}(t_{1},u_{1};k)u_{1}^{2i+1}du_{1}=0, \,\,\,\,\
i\in\{1,2\} \eqno(9.7)$$
 By (9.6) and (9.7) we obtain
$$(H_{k}g_{1})\left(t-\frac{1}{2}\right)=1+\int_{-\frac{1}{2}}^{\frac{1}{2}}u_{1}^{2i+1}du_{1}+
\int_{-\frac{1}{2}}^{\frac{1}{2}}K_{1}(t_{1},u_{1};k)u_{1}^{2i+1}du_{1}+
\int_{-\frac{1}{2}}^{\frac{1}{2}}K_{2}(t_{1},u_{1};k)u_{1}^{2i+1}du_{1}=$$
$$=1+\int_{-\frac{1}{2}}^{\frac{1}{2}}K_{2}(t_{1},u_{1};k)u_{1}^{2i+1}du_{1}=$$
$$=1+\left(\sqrt[k]{t_{1}^{3}+1}-1\right)\int_{-\frac{1}{2}}^{\frac{1}{2}}\psi_{3}(u_{1})u_{1}^{2i+1}du_{1}+
\left(\sqrt[k]{t_{1}^{5}+1}-1\right)\int_{-\frac{1}{2}}^{\frac{1}{2}}\psi_{4}(u_{1})u_{1}^{2i+1}du_{1}.$$
Case: $i=1$
$$(H_{k}g_{1})\left(t-\frac{1}{2}\right)=1+\left(\sqrt[k]{t_{1}^{3}+1}-1\right)\int_{-\frac{1}{2}}^{\frac{1}{2}}\psi_{3}(u_{1})u_{1}^{3}du_{1}+
\left(\sqrt[k]{t_{1}^{5}+1}-1\right)\int_{-\frac{1}{2}}^{\frac{1}{2}}\psi_{4}(u_{1})u_{1}^{3}du_{1}$$
By Lemma 9.2
$$\int_{-\frac{1}{2}}^{\frac{1}{2}}\psi_{3}(u_{1})u_{1}^{3}du_{1}=0, \,\,\
  \int_{-\frac{1}{2}}^{\frac{1}{2}}\psi_{4}(u_{1})u_{1}^{3}du_{1}=1.$$
  Then
  $$(H_{k}g_{1})\left(t-\frac{1}{2}\right)=1+\left(\sqrt[k]{1+t_{1}^{5}}-1\right)=
  \sqrt[k]{1+t_{1}^{5}}=\sqrt[k]{1+\left(t-\frac{1}{2}\right)^{5}}=g_{2}\left(t-\frac{1}{2}\right).$$
  Case: $i=2$
  $$(H_{k}g_{1})\left(t-\frac{1}{2}\right)=1+\left(\sqrt[k]{t_{1}^{3}+1}-1\right)\int_{-\frac{1}{2}}^{\frac{1}{2}}\psi_{3}(u_{1})u_{1}^{5}du_{1}+
\left(\sqrt[k]{t_{1}^{5}+1}-1\right)\int_{-\frac{1}{2}}^{\frac{1}{2}}\psi_{4}(u_{1})u_{1}^{5}du_{1}.$$
By Lemma 9.2 we get
$$(H_{k}g_{1})\left(t-\frac{1}{2}\right)=1+\left(\sqrt[k]{t_{1}^{3}+1}-1\right)=\sqrt[k]{t_{1}^{3}+1}=g_{1}\left(t-\frac{1}{2}\right).$$
Thus we have proved
$$(H_{k}g_{1})\left(t-\frac{1}{2}\right)=g_{2}\left(t-\frac{1}{2}\right), \,\,\,\
 (H_{k}g_{2})\left(t-\frac{1}{2}\right)=g_{1}\left(t-\frac{1}{2}\right).$$\\

 \textbf{Theorem 9.5} Let $k\geq k_{0}.$ The model
$$H(\sigma)=-\frac{1}{\beta}\sum_{<x,y>}ln\widetilde{K}\left(\sigma(x)-\frac{1}{2},\sigma(y)-\frac{1}{2};k\right)$$
 on the Cayley tree $\Gamma^{k}$ has at least four periodic Gibbs
 measures.\\

\textbf{References}\\

1.  Bleher, P.M.  and  Ganikhodjaev N.N.: On pure phases of the
Ising model on the Bethe lattice. {\it Theor. Probab. Appl.} {\bf
35} (1990), 216-227.

2.  Bleher, P.M.,  Ruiz, J.  and   Zagrebnov V.A.: On the purity
of the limiting Gibbs state for the Ising model on the Bethe
lattice. {\it Journ. Statist. Phys}. {\bf 79} (1995), 473-482.

3. Eshkabilov Yu.Kh, Haydarov F.H.,  Rozikov U.A.: Uniqueness of
Gibbs Measure for Models With Uncountable Set of Spin Values on a
Cayley Tree. {\it arXiv}:1202.1722v1 [math.FA]. To appear in {\it
Math. Phys. Anal. Geom}.

4. Eshkabilov Yu.Kh., Haydarov F.H., Rozikov U.A. : Non-uniqueness
of Gibbs Measure for Models With Uncountable Set of Spin Values on
a Cayley Tree {\it J.Stat.Phys}. {\bf 147} (2012), 779-794.

5. Ganikhodjaev, N.N. : On pure phases of the ferromagnet Potts
with three states on the Bethe lattice of order two. {\it Theor.
Math. Phys.} {\bf 85} (1990), 163--175.

6. Ganikhodjaev, N.N. and Rozikov, U.A. Description of periodic
extreme Gibbs measures   of some lattice model on  the Cayley
tree. {\it Theor. and Math. Phys}. {\bf 111} (1997), 480-486.

7. Ganikhodjaev, N.N. and Rozikov, U.A. : The Potts model with
countable set of spin values on a Cayley Tree. {\it Letters Math.
Phys.} {\bf 75} (2006), 99-109.

8. Ganikhodjaev, N.N. and Rozikov, U.A. On Ising model with four
competing interactions on Cayley tree. {\it Math. Phys. Anal.
Geom.} {\bf 12} (2009), 141-156.

9. H.O.Georgii, Gibbs measures and phase transitions, Second
edition. de Gruyter Studies in Mathematics, 9. Walter de Gruyter,
Berlin, 2011.

10. Krasnosel'ski, M.A.: Positive Solutions of Operator Equations.
Gos.Izd. Moskow (1969) (Russian).

11. Preston, C.: {\it Gibbs states on countable sets} (Cambridge
University Press, London 1974).

12. Rozikov, U.A. Partition structures of the Cayley  tree and
applications for describing periodic Gibbs distributions. {\it
Theor. and Math. Phys.} {\bf 112} (1997), 929-933.

13. Rozikov, U.A. and Eshkabilov, Yu.Kh.: On models with
uncountable set of spin values on a Cayley tree: Integral
equations. {\it Math. Phys. Anal. Geom.} {\bf 13} (2010), 275-286.

14. Sinai,Ya.G.: {\it Theory of phase transitions: Rigorous
Results} (Pergamon, Oxford, 1982).

15. Spitzer, F.: Markov random fields on an infinite tree, {\it
Ann. Prob.} {\bf 3} (1975), 387--398.

16. Suhov, Y.M. and  Rozikov, U.A.: A hard - core model on a
Cayley tree: an example of a loss network, {\it Queueing Syst.}
{\bf 46} (2004), 197--212.

17.  Zachary, S.: Countable state space Markov random fields and
Markov chains on trees, {\it Ann. Prob.} {\bf 11} (1983),
894--903.
\end{document}